\documentclass{aip-cp}

\usepackage[numbers]{natbib}
\usepackage{rotating}
\usepackage{graphicx}

\usepackage{savesym}
\savesymbol{iint}
\savesymbol{iiint}
\savesymbol{iiiint}
\savesymbol{idotsint}
\usepackage{amsmath}
\usepackage{amssymb}
\restoresymbol{AMS}{iint}
\restoresymbol{AMS}{iiint}
\restoresymbol{AMS}{iiiint}
\restoresymbol{AMS}{idotsint}
\usepackage{multirow}
\usepackage{textcomp}

\newtheorem{lemma}{Proposition}

\newcommand*{\QEDB}{\hfill\ensuremath{\square}}

\begin{document}

\title{Multistage Voting Model with Alternative Elimination}

\author[aff1]{Oleg A. Malafeyev}
\eaddress{malafeyevoa@mail.ru}
\author[aff1]{Denis Rylow\corref{cor1}}
\author[aff2,aff3]{Irina Zaitseva}
\eaddress{ziki@mail.ru}
\author[aff2]{Anna Ermakova}
\eaddress{dannar@list.ru}
\author[aff2]{Dmitry Shlaev}
\eaddress{shl-dmitrij@yandex.ru}

\affil[aff1]{St. Petersburg State University, 7/9 Universitetskaya nab., St. Petersburg, 199034, Russia.}
\affil[aff2]{Stavropol State Agrarian University, Zootekhnicheskiy lane 12, Stavropol, 355017, Russia.}
\affil[aff3]{Stavropol branch of the Moscow Pedagogical State University, Dovatortsev str. 66 g, Stavropol, 355042, Russia.}
\corresp[cor1]{Corresponding author: denisrylow@gmail.com}

\maketitle

\begin{abstract} 
The voting process is formalized as a multistage voting model with successive alternative elimination. A finite number of agents vote for one of the alternatives each round subject to their preferences. If the number of votes given to the alternative is less than a threshold, it gets eliminated from the game. A special subclass of repeated games that always stop after a finite number of stages is considered. Threshold updating rule is proposed. A computer simulation is used to illustrate two properties of these voting games.
\end{abstract}

\section{INTRODUCTION}

Voting is a procedure for consensus determination among group members with different opinions concerning alternatives. 
It plays an important role in modern societies: national assemblies and presidents are elected through voting process, chief executive officers of public companies are appointed by shareholders through a voting procedure. Voting itself can be viewed as an aggregating function of individual choices into collective compromise. The simplest forms of voting have been used since ancient times but the first formal analysis was given by Jean-Antoine the marquis de Condorcet and Jean-Charles de Borda. Both considered plurality voting and proposed different voting methods which became known as Condorcet method and Borda count \cite{condorcet}. Their works inspired further research and lead to such notable results as Dogson's method \cite{dogson}, Kemeny rule \cite{kem1,kem2}, Condorcet polling \cite{condorcet_polling}. 
Borda count and simple plurality voting are special cases of the voting model with scoring methods (see p.~133, \cite{moulin}). 
The widely discussed approval voting is a special case of grading voting or voting system with grading methods \cite{brams2,approval}. Grading voting is a nonranked voting procedure. Voters grade each alternative by assigning it an admissible number (usually from a finite set of consecutive natural numbers). Different alternatives may receive same grades from one agent. The winning alternative is selected subject to the grades and some rule. 




All models mentioned above are one stage models, the single alternative is chosen after only one cast of ballouts. Voting procedures that involve several rounds of voting are called multistage. Two round elections are multistage voting procedures and usually involve elimination of the least popular alternatives from the voting process before proceeding to the next round. A good overview of the literature is given in \cite{handbook}.

In this paper the voting process is formalized as a multistage voting model. Unlike runoff voting, alternatives are eliminated after each round of voting if they do not receive a sufficient number of votes, i.~e. do not reach a certain threshold. Agents vote again after each elimination of alternatives. The process stops when there is only one alternative left. Methods of the game theory are potent tools for modeling and analysis of socio-economic interaction \cite{Alferov2015,Neverova2015,Gureva2016115,Nepp2015,Wang,Malafeyev_scopus_2015_corrupt_model,Malafeyev_scopus_2014_electric_circutis_analogies_in_economy,Malafeyev_scopus_2014_postman_problem,Malafeyev_scopus_2015_stochastic_socio_economic_model,Malafeyev_scopus_2017_mean_field_game,Malafeyev_scopus_2016_parameter_mechanism_design,Malafeyev_Kolokoltsov_scopus_2010_understanding_game_theory_book,Malafeyev_scopus_2016_estimation_of_corruption_indicators,Malafeyev_scopus_2015_multicomponent_dynamics_single_sector_economy}. These methods can be applied to the voting process as well, especially when agents vote strategically. We formalize the voting process as a repeated voting game but do not perform game-theoretic analysis as the simplest case of agents voting sincerely according to their preferences is considered. 
The problems of stability and numerical analysis are presented in \cite{Kvitko3,Kvitko2,zyb1,zyb2,Smirnov,greece1,greece2}.

%
 The paper is organized as follows. The voting process is described in section ``Non-Formal Voting Process Description'' and formalized in the next section. The model is studied in section ``Analysis'': a subclass of games that stop after a finite number of stages is defined, a threshold updating rule is proposed. In section ``Computer Simulation'' the results of computer program are used to demonstrate two properties of the subclass of repeated games that stop after a finite number of stages as defined in the preceding sections. Conclusions are drawn in the last section.
 

\section{Non-Formal Voting Process Description}
Consider a voting process lasting several rounds with a finite number of alternatives that a finite number of voters participate in. Each voter has a finite number of votes which he/she assigns to one of the alternatives each round. Each alternative has a threshold --- a minimum number of votes necessary for the alternative to be passed to the next stage of the voting process. Once an alternative is eliminated, it can not be voted for again in any following rounds. The voting process stops when there is only one non-eliminated alternative left or less, i.~e. zero.
All voters are assumed to be rational and sincere. They vote according to their preferences. Preferences do no change during the voting process. 

\section{Formalization}
Consider a $n$-persons non-cooperative game $\gamma$.
Let $N$ be a set of all agents, where $N = \{ 1,2, \dots , n\}$. These agents are electors who vote for alternatives. Each agent possesses $v_i$ votes that he/she assigns to one of the alternatives. Let $V = (v_1, v_2, \dots, v_n)$. Without loss of generality $v_i \geq 1$, $v_i \in  \mathbb{N}$. Denote by $X = \{ x_1, x_2, \dots, x_m\}$ the set of all alternatives. Let $ \succ_i$ be a preference relation of the agent $i$. We say that the agent $i$ prefers alternative $x_s$ to the alternative $x_k$ iff $x_s \succ_i x_k$. Denote the minimum number of votes needed for the alternative $i$ to be passed to the next stage of the voting process by $f_i$. The value $f_i$ is a threshold. Then the game $\gamma$ is defined as a set $\gamma = (N, V, X, (\succ_i)_N, (f_i)_{i = 1}^{m})$. The strategy of the agent $i$ in the game $\gamma$ is denoted by $s_i \in X$. Let the profile in the game $\gamma$ be $s=(s_1, s_2, \dots, s_n) \in \underbrace{X \times \dots \times X }_{n \text{ times}}$. Let $r_i$ be a number of votes given to the alternative $i$ by the agents. Then we get $r_i = \sum_{j = 1}^{n} v_j \delta_i (s_j)$, where
$\delta_i (s_j)= 1$ if $s_j = x_i$ and $\delta_i (s_j)= 0$ if $ s_j \not= x_i $.
Without loss of generality assume that $f_i \leq \sum_{j = 1}^{n} v_j$ holds for at least some $i$. Otherwise the game has a trivial solution --- all alternatives are always eliminated.

Now we can construct a repeated voting game $\Gamma$. At each stage of the repeated voting game $\Gamma$ a game $\gamma$ is played by the agents $1,2,3, \dots, n$. Denote the number of stages of the game $\Gamma$ by $K$. Then we can write $\Gamma = (\gamma_k)_{k = 1}^{K}$, ${\gamma_k = (N, V, X^k, (\succ_i^k)_N, (f_i^k)_{i = 1}^{m_k}} )$. Alternative $i$ is eliminated from the voting process at the stage $k$ if $f^k_i > r_i^k$, where $ r_i^k = \sum_{j = 1}^{n} v_j    \delta_i (s_j^k) $, while 
$\delta_i (s_j^k)= 1$ if $  s_j^k = x^k_i$ and $\delta_i (s_j^k)= 0$ if $ s_j^k \not= x^k_i$.
Denote by $s^k =(s^k _1, s^k _2, \dots, s^k _n)$ a profile in the game $\gamma^k$. So $X^{k + 1} \subseteq X^{k }$.

\section{Analysis}

The formal definition of the repeated voting game given in the previous section does not guarantee that the number of stages $K$ is finite.
Consider a simple example.

\textbf{Example 1. } Let $N = \{ 1,2,3 \}$, $X^1 = \{ x^1_1, x^1_2, x^1_3\}$, $f_1^k = f_2^k = f_3^k = 2$, for any $k$, $V = (1,1,1)$. Assume that $x^1_1 \succ^k_1 x^1_2, x^1_2 \succ^k_1 x^1_3$, and $x^1_2 \succ^k_2 x^1_3, x^1_3 \succ^k_1 x^1_1$, and $x^1_3 \succ^k_3 x^1_2, x^1_2 \succ^k_3 x^1_1$ for any $k$. So the first, the second, and the third agent vote for $x^1_1$, $x^1_2$, $x^1_3$ respectively in accordance with their preferences. After the first round of voting no alternatives are eliminated as $r_1^1 = f_1^1$, $r_2^1 = f_2^1$, $r_3^1 = f_3^1$. This means that $\gamma^1= \gamma^2 =\dots=\gamma^K $ for any $K$. The game goes on ad infinitum.


\begin{lemma}
\label{one}
Let $\Gamma$ be a repeated voting game as formalized in section ``Formalization''. If $\sum_{i=1}^{m_k} f^k_i > \sum_{i=1}^{n} v_i$ for any $k$ 
 then 
\begin{enumerate}
\item at least one alternative is eliminated at every stage of the repeated voting game;
\item the length of the repeated voting game $K \leq m_1 - 1$.
\end{enumerate}
\end{lemma}

\emph{Proof. }
The second part of proposition immediately follows from the first part, as there are only $m$ alternatives and the game stops when there is only one alternative left or less. 

The proof of the first part is straightforward. Assume there exists such $k$ that no alternatives are eliminated from the game at the stage $k$. Then $r^k_i \geq f^k_i$, $i = 1, \dots, m_k$, where $r^k_i $ is the total number of votes given to the alternative $x_i^k$ by agents at the stage of game $k$. The sum $\sum_{i = 1}^{m_k } r^k_i = \sum_{i = 1}^{n } v_i$ because all agents vote for one of the alternatives. But then  $\sum_{i = 1}^{n } v_i = \sum_{i = 1}^{m_k }r_i \geq \sum_{i = 1}^{m_k }f_i$ which contradicts the conditions of the proposition.
\QEDB

Proposition \ref{one} allows us to consider a subclass of repeated voting games which stop after a finite number of stages regardless of the way agents vote. The proposition provides only sufficient conditions. It is possible to find other subclasses of repeated voting games that stop after a finite number of stages. Notice that the number of terms in $\sum_{i = 1}^{m_k } f_i^k$ decreases as $k$ increases while the sum $\sum_{i = 1}^{n } v_i$ stays constant. It means that in general case some terms of the sum $\sum_{i = 1}^{m_k } f_i^k$ have to be increased as the game progresses for $\sum_{i = 1}^{m_k }r_i^k < \sum_{i = 1}^{m_k }f_i^k$ to hold.
One way to solve this problem is to define a rule that will guarantee $\sum_{i = 1}^{m_k }r_i^k < \sum_{i = 1}^{m_k }f_i^k$ for every $m_k$. 

Let $a_{i}^{k} = r_i^k - f_i^k$ be the relative popularity of the alternative $x^k_i \in X^k$ among agents, $i = \overline{1, \dots, m_k}$ at the stage of the game $k$.  
Denote the previous lower index $i$ of the alternative $x^k_i$ by $I(i)$, so that $x^k_i = x^{k-1}_{I(i)}$, $f_i^k = f_{I(i)}^{k - 1}$

Now we can define the threshold updating rule to recursively find $f_i^k$ knowing $r_{I(i)}^{k-1}$, $a_{I(i)}^{k-1}$ and $f_{I(i)}^{k-1}$: 
$
f_i^k = f_{I(i)}^{k - 1} + \frac{a_{I(i)}^{k - 1} } { \sum_{j = 1}^{m_{k}} a_{I(j)}^{k - 1}} \Big( \sum_{j = 1}^{ m_{k -1}} f^{k - 1}_j  - \sum_{j = 1}^{m_{k}} f^{k - 1}_{I(j)}  \Big) .
$
This equation has clear interpretation. The threshold of each non-eliminated alternative $f_{I(i)}^{k - 1}$ is increased at every stage of the repeated voting game by a fraction of the sum of thresholds corresponding to the eliminated alternatives. The increases of $f_{I(i)}^{k - 1}$ are proportional to the relative popularity $a_{I(i)}^{k - 1}$ of the corresponding alternative $x_{I(i)}^{k - 1}$. Notice that the sum $\sum_{i=1}^{m_k} f_i^k$ stays constant: ${\sum_{i=1}^{m_j} f_i^j = \sum_{i=1}^{m_k} f_i^k}$, for any $ k, j \in \{ 1, \dots, K \}$. So, if $\sum_{j = 1}^{ m_{1}} f^{1}_j > \sum_{j = 1}^{n} v_j$. The repeated voting game always stops after a finite amount of stages according to the proposition \ref{one}.
Note however, that some games may stop after a finite number of stages with all alternatives being eliminated.

\begin{table}
\caption{Average game length}
\label{tab1}
\centering
\begin{tabular}{|c|c|c|c|c|c|c|c|c|c|c|c|c|c|c|c|c|}
 \hline
\multicolumn{11}{|c|}{agents} \\ \hline
\multirow{11}{*}{\rotatebox{90}{alternatives}} &~ & 2 & 4  & 8& 16 & 32 & 64 & 128 & 256 & 512  \\ \hline
& 10 &2.89 & 3.00 & 2.66 &2.15& 2.01& 2.00& 2.00& 2.00& 2.00 \\ 
& 20 &2.96 &3.00 &3.00 &2.99 &2.49 &2.20& 2.01 &2.00 &2.00 \\ 
& 40 &3.00 &3.00 &3.00 &3.00 &3.12 &2.90& 2.49 &2.16 &2.00 \\ 
& 80 &3.00 &3.00 &3.00 &3.00 &3.00 &3.37& 3.03 &3.05 &2.37 \\ 
& 160 &2.98 &3.00 &3.00 &3.00 &3.00 &3.00& 4.02 &3.09 &3.29  \\ 
& 320 &2.99 & 3.00& 3.00& 3.00& 3.00& 3.00 &3.00 &4.34& 3.10 \\ 
& 640 &3.00 & 3.00& 3.00& 3.00& 3.00& 3.00 &3.00 &3.00& 4.52 \\ 
& 1280 &3.00 & 3.00& 3.00& 3.00& 3.00& 3.00 &3.00 &3.00& 3.00 \\ 
& 2560&3.00 & 3.00& 3.00& 3.00& 3.00& 3.00 &3.00 &3.00& 3.00 \\ \hline
\end{tabular}
\end{table}

\section{Computer Simulation}
Consider a voting process with a fixed number of participating agents and a fixed number of alternatives. Assume that every agent has one vote. Take $f_i^1 = 2 \times \frac{ \text{number of agents} }{ \text{number of alternatives} }= \frac{n}{m_1}$, $i = \overline{1, m_1}$. All other $f_i^k$, $k = 2, 3, \dots$, are recursively found. 
Then it can be demonstrated that for any distribution of agents' preferences:
\begin{enumerate}
\item increments in the number of participating agents at first increase and then decrease the length of games on average if the initial number of agents is small;
\item decrements in the number of available alternatives increase and then decrease the length of games on average if the number of agents is sufficiently large.
\end{enumerate}
To illustrate these two points we provide the results of computer simulation. One hundred experiments are performed and the average length of games is computed for different numbers of alternatives and different numbers of voting agents. Agent's preferences are uniformly distributed. The results are listed in the table \ref{tab1}; rows correspond to the number of alternatives, columns correspond to the number of voting agents.


\section{Conclusion}
Voting process is formalized as a $n$-persons non-cooperative multistage repeated voting game. 
Sufficient conditions that allow us to consider repeated voting games that stop after a finite number of stages are given.
Threshold updating rule is formulated and applied recursively to the threshold of every alternative at every stage of the repeated voting game except the first stage. A computer simulation is used to illustrate two properties of repeated voting games that satisfy sufficient conditions.





\hyphenpenalty=10000
\nocite{*}
\bibliography{malafeyev_scopus_web_of_science(Science_Index)}%
\bibliographystyle{ugost2008}%

\end{document}